\documentclass[12pt]{article}
\usepackage[all]{xy}
\usepackage{latexsym,amsmath,amssymb,amsthm
}
\newtheorem{thm}{Theorem}[section]

\newtheorem{lem}[thm]{Lemma}
\newtheorem{cor}[thm]{Corollary}
\newtheorem{rem}[thm]{Remark}

\theoremstyle{definition}  
\newtheorem{defn}[thm]{Definition}
\newcommand{\arrow}{{\longrightarrow}}

\def\id{{\rm id}}

\def\Pos{{\rm Pos}}
\def\Con{{\rm Con}}
\date{}
\title{\bf Characterization of Pomonoids by Properties of Generators }
\author{{\bf Setareh Irannezhad}\and {\bf Ali Madanshekaf } \\
Department of Mathematics\\Faculty of Mathematics, Statistics and Computer Science\\
Semnan University\\
Semnan\\
Iran\\ emails: amadanshekaf@semnan.ac.ir\\
setareh.irannezhad@students.semnan.ac.ir }
\date{ November 8, 2014}
\begin{document}
\maketitle
\begin{abstract}
The study of flatness properties of ordered monoids acting on posets
was initiated by S.M. Fakhruddin in the 1980's. Although there exist
many papers which investigate various properties of $S$-posets
(posets equipped with a compatible right action of an ordered monoid
$S$)  from free to torsion free, among them generators, there seems
to be known very little. In 2008, Laan characterized generators in
the category {\bf Pos}-$S$ of all $S$-posets with monotone
action-preserving maps between them. His characterization is similar
to the case of acts over monoids. We attempt here to collect the
knowledge on generators in the category {\bf Pos}-$S$ and to apply
this to proceed on the questions of homological classification of
ordered monoids, that is results of the type: all generators in the
category  {\bf Pos}-$S$, satisfy a flatness property if and only if
$S$ has a certain property.
\end{abstract}
{\bf AMS} {\it {\bf Subject classification}}: 20M50; 06F05. \\
{\it {\bf Key words} }: Ordered monoid, $S$-poset, Flatness, Generator
\section{Introduction and Preliminaries}
Over the past several decades, a large body of literature has accumulated that investigates
what are usually called flatness properties of acts over monoids. This research has to do with such notions as
flatness and projectivity, and is usually connected with preservation properties of the functor
$A_S\otimes {}_{S} -$ (from the category of left $S$-acts into the category of sets),
for a right act $A_S$ over a monoid $S$. The monograph by Kilp et al. (2000)
furnishes an exhaustive and reasonably current treatment of this area.

In the 1980's, Fakhruddin published two articles (1986, 1988)
devoted to tensor products and flatness properties in the context of
ordered monoids acting (monotonically in both arguments) on ordered
sets (that is, $S$-posets). Further investigations along these lines
have lain dormant until the recent appearance of the articles Shi et
al. (2005), Bulman-Fleming and Laan (2005), Shi (2005), and
Bulman-Fleming and Mahmoudi (2005).

Kilp and Knauer, in~\cite{character}, investigated monoids over
which all generators in the category of right $S$-acts satisfy a
special flatness property. Continuing this study, Sedaghatjoo,
in~\cite{sed:gen}, investigated the corresponding problem for weak
flatness, Condition $(E)$, and regularity in the category of right
$S$-acts. The purpose of the present article is to discuss pomonoids
over which all generators in the category {\bf Pos}-$S$ satisfy a
special property, that is related to flatness like Condition $(E)$,
weak po-flatness, principal weak po-flatness, po-torsion freeness
and $I$-regularity.

A monoid $S$ that is also a partially ordered set, in which the binary operation and the order relation are compatible,
is called a \emph{partially ordered monoid}, or shortly \emph{pomonoid}. A \emph{right $S$-poset}, often denoted $A_S$,
is a non-empty set $A$ equipped with a partial order $\leq$ and a right action
$A\times S\arrow A$, $(a, s)\rightsquigarrow as$, that satisfies:
\begin{itemize}
    \item[(1)] $a(ss')=(as)s'$,
    \item[(2)]  $a1=a$,
    \item[(3)] $a\leq a'$ implies $as\leq a's$,
    \item[(4)] $s\leq s'$ implies $as\leq as'$
\end{itemize}
for all $a, a'\in A$ and $s, s'\in S$. Left $S$-posets ${}_SB$ are
defined analogously and the left or right $S$-posets form
categories, denoted by $S$-{\bf Pos} and {\bf Pos}-$S$
respectively, in which the morphisms are the functions that preserve
both the action and the order which is called
\emph{$S$-homomorphism}. Let $A_S$ and $B_S$ be two $S$-posets, the
set of all $S$-homomorphisms from $A_S$ into $B_S$ will be denoted
by $\Pos_S(A_S,B_S)$. An \emph{$S$-subposet} of a right $S$-poset
$A_S$ is a subset of $A_S$ that is closed under the $S$-action. The
definition of ideal is the same as for the act case, that is a
\emph{right ideal} of a partially ordered monoid $S$ is a nonempty
subset $I$ of $S$ such that $IS\subseteq I$. If $X$ is a subset of a
poset $P$, $(X]:=\{p\in P~|~p\leq x, \text{for some}~ x\in X\}$ is
the down-set of $X$.

In Bulman-Fleming and Laan (2005) a careful treatment was given of factor $S$-posets, the essentials of which
we now repeat. A \emph{congruence} on an $S$-poset $A_S$ is an $S$-act congruence $\theta$
such that has the further property that the factor act $A_S/\theta$ can be equipped with a compatible order
so that the natural map $A_S\arrow A_S/\theta$ is an $S$-homomorphism. An alternative description is as follows:
if $\sigma$ is any binary relation on $A_S$ we write $a\leq_{\sigma}a'$ if a so-called \emph{$\sigma$-chain}
\begin{equation}\label{}
 a\leq a_1 ~\sigma~ a'_1\leq a_2~\sigma~ a'_2\leq \cdots ~\leq a_m\sigma~ a'_m\leq a'
\end{equation}
from $a$ to $a'$ exists in $A_S$. It can be shown that an $S$-act
congruence $\theta$ on an $S$-poset $A_S$ is an $S$-poset congruence
precisely when $a\leq_{\theta}a'\leq_{\theta}a$ implies $a~\theta~
a'$ for all $a, a'\in A_S$. In speaking of an $S$-poset congruence
$\theta$, we will always  specify the corresponding order on
$A_S/\theta$ that we have in mind. (A given $S$-act congruence can
of course give rise to many different factor $S$-posets: for
example, if $A_S$ is any $S$-poset with compatible order $\leq$ then
the identity relation $\Delta$ is an $S$-poset congruence on $A_S$
and the corresponding factor acts are simply the S-posets $(A_S,
\leq ')$ where the relation $\leq '$ is compatible and finer than
the relation $\leq$).

Let $A_S$ be an $S$-poset and $\alpha$ is a binary relation on $A_S$
that is reflexive, transitive, and compatible with the $S$-action.
Here is a  description of the smallest $S$-poset congruence $\theta$
on $A_S$ that contains $\alpha$, together with a natural
accompanying order on $A_S/\theta$. One may check that the relation
$\theta$ defined on $A_S$ by
\begin{equation}\label{}
  a~\theta~a' ~\text{if and only if}~ a\leq_{\alpha} a'\leq_{\alpha} a
\end{equation}
is an $S$-act congruence on $A_S$, a suitable order relation on
$A_S/\theta$ being
\begin{equation}\label{}
 [a]_{\theta}\leq [a']_{\theta}~ \text{if and only if}~ a\leq_{\alpha} a'.
\end{equation}
Furthermore, if $\eta$ is any $S$-poset congruence on $A_S$ such
that $\alpha \subseteq \eta$, then $\theta \subseteq \eta$ as well.
The relation $\theta$ is called the \emph{$S$-poset congruence
generated by $\alpha$}. In particular, if $H\subseteq A_S\times A_S$
and $\alpha$ is the $S$-act congruence on $A_S$ generated by $H$,
the corresponding $S$-poset congruence $\theta$ will be denoted
$\theta(H)$ and will be called the \emph{$S$-poset congruence
generated by $H$}.

The next idea that is given in~\cite{lazard}, is even more important because it is used for
the definition of monocyclic $S$-poset congruences in this paper, furthermore,
recall that this idea had an important role in the proof of Stenstrom-Govorov-Lazard Theorem in~\cite{lazard}.\\
Let $A_S$ be an $S$-poset and let $H\subseteq A_S\times A_S$. Then there exists an $S$-poset congruence
$\nu(H)$ on $A_S$ having the properties
\begin{itemize}
    \item[(1)] if $(h, h')\in H$, then $[h]_{\nu(H)}\leq [h']_{\nu(H)}$; and
   \item[(2)] if $\rho$ is any $S$-poset congruence on $A_S$ such that $[h]_{\rho}\leq [h']_{\rho}$
   whenever $(h, h')\in H$, then $\nu(H)\subseteq \rho$.
\end{itemize}
The relation $\nu(H)$ is called the \emph{$S$-poset congruence on $A_S$ induced by $H$}.
Define a relation $\alpha(H)$ on $A_S$ by $a~\alpha(H)~ a'$ if and only if $a=a'$ or
\begin{equation}\label{}
   \begin{array}{rcl}
a&=&h_1s_1\\
h'_1s_1&=&h_2s_2\\
&\cdots&\\
h'_ns_n&=&a'.
\end{array}
\end{equation}
for some $(h_i, h'_i)\in H$ and $s_i\in S$. Note that the relation
$\alpha(H)$ is transitive, reflexive and $a~\alpha(H)~ a'$ implies
$as~\alpha(H)~ a's$ for every $a, a'\in A_S$ and $s\in S$.
Therefore, the relation $\nu(H)$ defined on $A_S$ by
$$a~\nu (H) ~a' ~\text{if and only if}~ a\leq _{\alpha (H)} a'\leq _{\alpha (H)} a$$
is the smallest $S$-poset congruence containing $\alpha(H)$, where $[a]_{\nu (H)}\leq [a']_{\nu (H)}$
exactly when $a\leq_{\alpha (H)} a'$. We note in particular that $\theta (H)=\nu (H\cup H^{op})$ for any
$H\subseteq A_S\times A_S$.

An $S$-poset $A_S$ is called \emph{(principally) weakly po-flat} if
the functor $A_S\otimes {}_{S} -$ preserves embeddings of
(principal) left ideals into $S$ (Note that, for $S$-posets and for
posets, the monomorphisms are the injective, monotonic morphisms,
whereas the embeddings are the order-embeddings).

Conditions $(P)$ and $(E)$ for $S$-posets are formulated as follows (see for example Bulman-Fleming and Laan, 2005).\\
An $S$-poset $A_S$ satisfies \emph{Condition $(P)$} if, for all $a,
a'\in A_S$ and $s, t\in S$, $as\leq a't$ implies $a=a''u$, $a'=a''v$
for some $a''\in A_S$, $u, v\in S$ with $us\leq vt$. The $S$-poset
$A_S$ satisfies \emph{Condition $(E)$} if, for all $a\in A_S$ and
$s, t\in S$, $as\leq at$ implies $a=a'u$ for some $a'\in A_S$, $u\in
S$ with $us\leq ut$. An $S$-poset $B_S$ is said to be \emph{strongly
flat} if it satisfies both Conditions $(E)$ and $(P)$.

An object $G_S$ in the category {\bf Pos}-$S$ is called a
\emph{generator} if the functor $\Pos_S(G_S, -)$ is faithful, i.e.
for any $S$-posets $X_S$, $Y_S$ and any $f, g \in \Pos_S(X_S, Y_S)$
with $f\neq g$ there exists $\alpha\in \Pos_S(G_S, X_S)$ such that
$f\alpha\neq g\alpha$.  In~\cite{gen} generators have been
characterized in the category {\bf Pos}-$S$ by Laan. Recalling that
a monoid $S$ is \emph{right reversible} if $Ss\cap St\neq \emptyset$
for all $s, t\in S.$ We shall call a pomonoid $S$ \emph{weakly right
reversible} in case $Ss\cap (St]\neq \emptyset$ for all $s, t\in S$.

The reader is referred to Howie (1995) and Kilp et al. (2000), respectively,
for information on general semigroup theory and flatness properties of $S$-acts that is not fully explained here.

Using~\cite[Theorem 2.1]{gen}, it is easy to see  that $S_S\times
A_S$ is a generator, for every $S$-poset $A_S$. So we have the next
fundamental theorem which yields characterizations in some theorems
in sequel and its proof is similar to the corresponding proof for
$S$-acts.
\begin{thm}\label{gen}
Let $S$ be a pomonoid and $\mathcal{P}$ be an $S$-poset property which is preserved under retractions.
The following
assertions are equivalent:
\begin{itemize}
      \item[(i)] all generators satisfy property $\mathcal{P}$;
      \item[(ii)] $S_S\times A_S$ satisfies property $\mathcal{P}$ for every right $S$-poset $A_S$;
      \item[(iii)] a right $S$-poset $A_S$ satisfies property $\mathcal{P}$ if $\Pos_S(A_S,S_S)\neq \emptyset$.
\end{itemize}
\end{thm}
\vspace{.2cm} \noindent {\bf Proof.}
$(i) \Rightarrow (ii)$ Suppose that all generators satisfy property $\mathcal{P}$. Let $A_S$ be a right $S$-poset.
Since $S_S\times A_S$ is a generator, using our assumption $S_S\times A_S$ satisfies property $\mathcal{P}$.

$(ii) \Rightarrow (iii)$ Let $A_S$ be a right $S$-poset such that
$\Pos_S(A_S,S_S)\neq \emptyset$. Take $f\in \Pos_S(A_S,S_S),$ then
by the universal property of product, there exists $h:
A_S\rightarrow S_S \times A_S$ making diagram
\begin{displaymath}
\xymatrix{
&  A_S  \ar[dl]_f   \ar@{-->}[d]^h \ar[dr]^{\id_A}   & \\
S_S  &          {S_S \times A_S}  \ar[l]_{\pi_1} \ar[r]^{\pi_2}           & A_S}
\end{displaymath}
commutative. Thus, $A_S$ is a retract of $S_S\times A_S$ which by
assumption satisfies property $\mathcal{P}$.

$(iii) \Rightarrow (i)$ Let $A_S$ be a generator. By~\cite[Theorem
2.1]{gen}, we conclude $\Pos_S(A_S,S_S)\neq \emptyset$, so $A_S$
satisfies property $\mathcal{P}$. $\blacksquare$

If $S$ is a pomonoid, the cartesian product $S^I$ is a right and left $S$-poset equipped with the order
and the action componentwise where $I$ is a nonempty set. Moreover, $(s_i)_{i\in I}\in S^I$
is dented simply by $(s_i)$, the right $S$-poset $S\times S$ will be denoted by $D(S)$ and will be called
\emph{diagonal $S$-poset} over $S$. By the fact that $\Pos_S(S^I, S)\neq \emptyset$ for each nonempty set $I$,
the following corollary is obtained.
\begin{cor}\label{S^I}
Let $S$ be a pomonoid and $\mathcal{P}$ be an $S$-poset property which is preserved under retractions.
If all generators satisfy property $\mathcal{P}$ then for each nonempty set $I$, $S^I$ satisfies property $\mathcal{P}$.
\end{cor}

Since for each right ideal $K$ of $S$, $\Pos_S(K_S, S_S)\neq \emptyset$, by the former theorem the following
corollary is deduced.
\begin{cor}
Let $S$ be a pomonoid and $\mathcal{P}$ be an $S$-poset property which is preserved under retractions.
If all generators satisfy property $\mathcal{P}$ then all right ideals of $S$ satisfy property $\mathcal{P}$.
\end{cor}
\section{Pomonoids over which  all generators  satisfy \\ Condition $(E)$}
In this section similar to the $S$-act case, we investigate that when all generators satisfy Condition $(E)$.
To reach this goal it will be useful to show that for all monocyclic $S$-posets the
concepts of strong flatness and projectivity are the same.\\
By analogy with the unordered case a right congruence $\rho$ on pomonoid $S$ is said to be \emph{monocyclic}
if it is induced by a single pair of elements $(s, t)$, $s, t\in S$ and is then denoted by $\nu(s, t)$.
A right factor $S$-poset of a pomonoid $S$ by a monocyclic right congruence is called a
\emph{monocyclic right $S$-poset} (see Definition 1.4.18 in~\cite{K:act} for the unordered case).
\begin{lem}\label{induced}
Let $A_S$ be an $S$-poset and $H\subseteq A_S\times A_S$. Then $a~\nu(H)~b$ if and only if $a=b$, or there exist
$s_1, s_2, \cdots , s_n, t_1, t_2, \cdots , t_m\in S$ such that
\begin{align}\label{nu(H)}
a\leq c_1s_1,~d_1s_1\leq c_2s_2,~d_2s_2\leq c_3s_3,~\cdots,~d_ns_n\leq b;\nonumber \\
b\leq p_1t_1~q_1t_1\leq p_2t_2,~q_2t_2\leq p_3t_3,~\cdots,~q_mt_m\leq a,
\end{align}
where $(c_i, d_i)\in H$, $(p_j, q_j)\in H$ for $i=1, 2, \cdots, n$, $j=1, 2, \cdots, m$, (\emph{$\nu(H)$
is the $S$-poset congruence on $A_S$ induced by $H$}).
\end{lem}
\vspace{.2cm} \noindent {\bf Proof.} We define a relation $\sigma$
on $A_S$ as follows: for any $a, b\in A_S, a~\sigma~b$ if and only
if $a=b$, or a system of inequalities~\eqref{nu(H)} exists. It is
easy to check that the relation $\sigma$ is an equivalence on $A_S$
which is compatible with the $S$-action. To see that $\sigma$ is an
order congruence, assume $a, a_i, a'_i\in A_S$ for $1\leq i\leq n$,
and suppose
$$a\leq a_1~\sigma ~a'_1\leq a_2~\sigma~a'_2\leq \cdots \leq a_j~\sigma~a'_j\leq \cdots
\leq a_{n-1}~\sigma~a'_{n-1}\leq a_n~\sigma~a'_n\leq a$$
then, for each $k\in \{1, \cdots, n\}$ we have a system of inequalities
\begin{equation*}
\begin{matrix}
a_k\leq c_{k, 1}s_{k, 1},\\
d_{k, 1}s_{k, 1}\leq c_{k, 2}s_{k, 2},\\
d_{k, 2}s_{k, 2}\leq c_{k, 3}s_{k, 3},\\
\cdots\\
d_{k, m_k}s_{k, m_k}\leq a'_k;
\end{matrix} \qquad
\begin{matrix}
a'_k\leq p_{k, 1}t_{k, 1},\\
q_{k, 1}t_{k, 1}\leq p_{k, 2}t_{k, 2},\\
q_{k, 2}t_{k, 2}\leq p_{k, 3}t_{k, 3},\\
\cdots\\
q_{k, l_k}t_{k, l_k}\leq a_k,
\end{matrix}
\end{equation*}
where $(c_{i, j}, d_{i, j})\in H$, $(p_{i, j}, q_{i, j})\in H$. Hence
\begin{equation*}
\begin{matrix}
a_1\leq c_{1, 1}s_{1, 1},\\
d_{1, 1}s_{1, 1}\leq c_{1, 2}s_{1, 2},\\
d_{1, 2}s_{1, 2}\leq c_{1, 3}s_{1, 3},\\
\cdots\\
d_{1, m_1}s_{1, m_1}\leq a'_1;
\end{matrix} \qquad
\begin{matrix}
a_2\leq c_{2, 1}s_{2, 1},\\
d_{2, 1}s_{2, 1}\leq c_{2, 2}s_{2, 2},\\
d_{2, 2}s_{2, 2}\leq c_{2, 3}s_{2, 3},\\
\cdots\\
d_{k, m_k}s_{k, m_k}\leq a'_k;
\end{matrix}
\end{equation*}
\begin{equation*}
\begin{matrix}
a_{k+1}\leq c_{k+1, 1}s_{k+1, 1},\\
d_{k+1, 1}s_{k+1, 1}\leq c_{k+1, 2}s_{k+1, 2},\\
d_{k+1, 2}s_{k+1, 2}\leq c_{k+1, 3}s_{k+1, 3},\\
\cdots\\
d_{k+1, m_{k+1}}s_{k+1, m_{k+1}}\leq a'_{k+1};
\end{matrix} \qquad
\begin{matrix}
a_{k+2}\leq c_{k+2, 1}s_{k+2, 1},\\
d_{k+2, 1}s_{k+2, 1}\leq c_{k+2, 2}s_{k+2, 2},\\
d_{k+2, 2}s_{k+2, 2}\leq c_{k+2, 3}s_{k+2, 3},\\
\cdots\\
d_{n, m_n}s_{n, m_n}\leq a'_n.
\end{matrix}
\end{equation*}
On the other hand, $a\leq a_1, a'_1\leq a_2, \cdots, a'_k\leq a_{k+1}, \cdots, a'_n\leq a$.
Therefore $a~\sigma ~a'_k$ holds for each $k$. Because $a_k~\sigma~a'_k$, we obtain $a~\sigma~a_k$.
Thus every closed $\sigma$-chain is contained in a single equivalence class of $\sigma$, so
$\sigma$ is an order congruence on $A_S$ by~\cite[Theorem 1.1]{Fakhruddin1}. Suppose $a~\nu(H)~b$,
then $a\leq_{\alpha(H)}b$, $b\leq_{\alpha(H)}a$. This means that there are the following $\alpha(H)$-chains
$$a\leq x_1~\alpha(H)~x'_1\leq x_2~\alpha(H)x'_2\leq \cdots \leq x_z~\alpha(H)~x'_z\leq b$$
$$b\leq y_1~\alpha(H)~y'_1\leq y_2~\alpha(H)y'_2\leq \cdots \leq y_t~\alpha(H)~y'_t\leq a.$$
Then for each $l\in \{1, 2, \cdots, z\}$, we have a system of equalities
\begin{equation*}
\begin{matrix}
x_l=h_{l, 1}s_{l, 1},\\
h'_{l, 1}s_{l, 1}=h_{l, 2}s_{l, 2},\\
\cdots\\
h'_{l, n_l}s_{l, n_l}=x'_l,
\end{matrix}
\end{equation*}
where $(h_{i, j}, h'_{i, j})\in H$. Hence,
\begin{align}
a\leq x_1=h_{1, 1}s_{1, 1},~~h'_{1, 1}s_{1, 1}\leq h_{1,2}s_{1, 2}, ~\cdots~, h'_{1, n_1}s_{1, n_1}=x'_1\leq\nonumber\\
x_2=h_{2, 1}s_{2, 1},~~h'_{2, 1}s_{2, 1}\leq h_{2, 2}s_{2, 2},~
\cdots~, h'_{2, n_2}s_{2, n_2}=x'_2\leq
\cdots x'_{l-1}\leq \nonumber\\
x_l=h_{l, 1}s_{l, 1},~~h'_{l, 1}s_{l, 1}\leq h_{l, 2}s_{l, 2},~
\cdots~, h'_{l, n_l}s_{l, n_l}=x'_l \leq x_{l+1}
\cdots \leq \nonumber\\
x_z=h_{z, 1}s_{z, 1},~~h'_{z, 1}s_{z, 1}\leq h_{z, 2}s_{z, 2,}~
\cdots~, h'_{z, n_z}s_{z, n_z}=x'_z\leq b. \nonumber
\end{align}
Similarly, there is  such a chain from $b$ to $a$, by considering
$\alpha(H)$-chain from $b$ to $a$. So $a~\sigma~ b.$ Consequently
$\nu(H)\subseteq \sigma$. On the other hand, if $a~\sigma~b$, then
there exists a system of inequalities~(\ref{nu(H)}). Therefore
$a\leq c_1s_1~\alpha(H)~d_1s_1\leq c_2s_2~\alpha(H)~d_2s_2 \leq
c_3s_3~\alpha(H)~\cdots~\alpha(H)~d_ns_n\leq b$. Thus
$a\leq_{\alpha(H)} b$. Similarly $b\leq_{\alpha(H)} a$. By the
definition of $\nu(H)$, $a~\nu(H)~b$. This means that $\sigma
\subseteq \nu(H)$. Thus $\sigma =\nu(H)$.
\begin{rem}\label{order}
The order relation on $A_S/\nu(H)$ can be defined as follows: $[a]_{\nu(H)}\leq [b]_{\nu(H)}$ if and
only if $a\leq b$, or there exist $s_1, s_2, \cdots, s_n\in S_S$ such that
$$a\leq c_1s_1,~d_1s_1\leq c_2s_2,~d_2s_2\leq c_3s_3,~\cdots,~d_ns_n\leq b$$
where $(c_i, d_i)\in H$ for $i=1, 2, \cdots, n$.
\end{rem}
\begin{thm}\label{monocyclic}
Every strongly flat monocyclic $S$-poset is projective.
\end{thm}
\vspace{.2cm} \noindent {\bf Proof.}
Let $s, t\in S_S$ and let $\nu (s, t)$ be $S$-poset congruence induced by $H=\{(s, t)\}$.
Assume the right $S$-poset $S_S/\nu (s, t)$ is strongly flat. Since for cyclic $S$-posets strong flatness
and Condition $(E)$ coincide (see ~\cite[Lemma 2.6]{sflat poflat}),
it follows from ~\cite[Proposition 3.12]{cyclic} that there exists $u\in S$ such that $us\leq ut$ and $u~\nu (s, t) ~1$.
Then by Lemma~\ref{induced}, $u=1$ or there exist $s_1, s_2, \cdots, s_n, t_1, t_2, \cdots, t_m\in S_S$ such that
\begin{center}
$u\leq c_1s_1,~d_1s_1\leq c_2s_2,~d_2s_2\leq c_3s_3,~\cdots,~d_ns_n\leq 1$;\\
$1\leq p_1t_1,~q_1t_1\leq p_2t_2,~q_2t_2\leq p_3t_3,~\cdots,~q_mt_m\leq u$,
\end{center}
where $(c_i, d_i)\in \{(s, t)\}$, $(p_j, q_j)\in \{(s, t)\}$ for $i=1, 2, \cdots, n$, $j=1, 2, \cdots, m$.
Since $uc_i\leq ud_i$ and $up_j\leq uq_j$ for $i=1, 2, \cdots, n$, $j=1, 2, \cdots, m$, we have
\begin{center}
$u^2\leq uc_1s_1\leq ud_1s_1\leq uc_2s_2\leq~\cdots~\leq ud_ns_n\leq u\leq up_1t_1\leq uq_1t_1
\leq up_2t_2\leq ~\cdots~\leq uq_mt_m\leq u^2$.
\end{center}
So $u^2=u$, thus $u$ is an idempotent. Next, let $x, y\in S_S$ and $[x]_{\nu (s, t)}\leq [y]_{\nu (s, t)}$.
Then, by Remark~\ref{order}, $x\leq y$ or there exist $s_1, s_2, \cdots, s_n\in S_S$ such that
\begin{center}
$x\leq c_1s_1,~d_1s_1\leq c_2s_2,~d_2s_2\leq c_3s_3,~\cdots,~d_ns_n\leq y$
\end{center}
where $(c_i, d_i)\in \{(s, t)\}$ for $i=1, 2, \cdots, n$. Now
\begin{center}
$ux\leq uc_1s_1\leq ud_1s_1\leq uc_2s_2\leq ud_2s_2\leq uc_3s_3\leq \cdots \leq ud_ns_n\leq uy$.
\end{center}
Thus, the right $S$-poset $S_S/\rho$ is projective
by~\cite[Proposition 3.10]{cyclic}. $\blacksquare$

A right congruence $\rho$ on $S_S$ is called a \emph{right subannihilator congruence} if $\rho\leq \ker \lambda_s$
(the left translation by $s$) for some $s\in S$ (see also~\cite{sed:gen} for unordered case).
It is easy to check that a right congruence $\rho$ on $S_S$ is a right subannihilator
congruence if and only if $\Pos_S(S_S/\rho, S_S)\neq \emptyset$.
\begin{defn}
Let $S$ be a pomonoid and $a, b\in S$. We define \emph{$l(a, b)=\{s\in S~|~ sa\leq sb\}$}.
\end{defn}

By the above definition, $l(a, b)=\emptyset$ or it is a left ideal
of $S$.
\begin{lem}
Let $S$ be a pomonoid. Then for each pair $(a, b)\in S\times S$ such that $l(a, b)\neq \emptyset$, the
monocyclic right congruence on $S_S$ induced by the pair $(a, b)$ is a right subannihilator congruence.
\end{lem}
\vspace{.2cm} \noindent {\bf Proof.} Suppose $(a, b)\in S\times S$
and $l(a, b)\neq \emptyset$ then there exists $s\in S$ such that
$sa\leq sb$. So $(a, b)\in K_{\lambda_s}$ where $K_{\lambda_s}=\{(x,
y)\in S\times S~|~ sx\leq sy\}$ is the directed kernel of
$\lambda_s$. By the~\cite[Proposition 2.3]{lazard}, all
$K_{\lambda_s}$, $\alpha(K_{\lambda_s})$ and
$\leq_{\alpha(K_{\lambda_s})}$ coinside. Therefore
$a\leq_{\alpha(K_{\lambda_s})} b$ so $[a]_{\ker \lambda_s}\leq
[b]_{\ker \lambda_s}$, since $\nu (K_{\lambda_s})=\ker \lambda_s$.
But $\nu (a, b)$ is the smallest $S$-poset congruence with this
property, so $\nu (a, b)\subseteq \ker \lambda_s$. $\blacksquare$

The following lemma is an easy consequence of definitions.
\begin{lem}\label{retE}
Let $S$ be a pomonoid. Then any retract of an $S$-poset satisfying Condition $(E)$ satisfies Condition $(E)$.
\end{lem}
Now we have the main result of this section.
\begin{thm}\label{mainE}
For a pomonoid $S$ the following assertions are equivalent:
\begin{itemize}
      \item[(i)] all generators satisfy Condition $(E)$;
      \item[(ii)] $S_S\times A_S$ satisfies Condition $(E)$ for each right $S$-poset $A_S$;
      \item[(iii)] a right $S$-poset $A_S$ satisfies Condition $(E)$ if $\Pos_S(A_S, S_S)\neq \emptyset$;
      \item[(iv)] for each right subannihilator congruence $\rho$, $S_S/\rho$ satisfies Condition $(E)$;
      \item[(v)] for each $x,y\in S$ the monocyclic right $S$-poset $S_S/\nu (x, y)$ satisfies Condition $(E)$
in case $l(x, y)\neq \emptyset$;
     \item[(vi)] for each $x,y\in S$ the monocyclic right $S$-poset $S_S/\nu (x, y)$ is projective or equivalently,
$\nu (x, y)=\ker \lambda_e$ for an idempotent $e\in S$ in case $l(x, y)\neq \emptyset$.
\end{itemize}
\end{thm}
\vspace{.2cm} \noindent {\bf Proof.} Assertions $(i)$, $(ii)$ and
$(iii)$ are equivalent by Theorem~\ref{gen} and Lemma~\ref{retE}.
Implications $(iii) \Rightarrow (iv)$ and $(iv) \Rightarrow (v)$ is
similar to the act case (see~\cite{sed:gen}).

$(v) \Rightarrow (vi)$ Suppose $x, y\in S$ and $l(x, y)\neq \emptyset$.
Then by the assumption $S_S/\nu (x, y)$ satisfies Condition $(E)$.
It is proved in Theorem~\ref{monocyclic} that every monocyclic right $S$-poset satisfying Condition
$(E)$ is projective. So $S_S/\nu (x, y)$ is projective or equivalently,
$\nu (x, y)=\ker \lambda_e$ for an idempotent $e\in S$ (see~\cite[Proposition 3.2]{Inde}).

$(vi) \Rightarrow (i)$ Let $A_S$ be a generator and $f: A_S\arrow S_S$ be an epimorphism.
Let $ax\leq ay$ for $a\in A_S$, $x, y\in S$ and $f(a)=s$. Then $f(a)x\leq f(a)y$ or $sx\leq sy$ so
$s\in l(x, y)$ and by our assumption $S/\nu (x, y)$ is projective.
Therefore there exists an idempotent $e\in S$ such that $\nu(x, y)=\ker \lambda_e$.
Since $ax\leq ay$, then $(x, y)\in K_{\lambda_a}$. So $\ker \lambda_e=\nu (x, y) \leq \ker \lambda_a$.
Since $(1, e)\in \nu (x, y)$ then $(1, e)\in \ker \lambda_a$ and we
get $a=ae$. By the definition of $\nu(x, y)$  the $S$-poset
congruence induced by a single pair $(x, y)$, we have $[x]_{\nu(x,
y)}\leq [y]_{\nu(x, y)}$, since $S_S/\nu (x, y)\cong eS$ then
$ex\leq ey$. Thus $A_S$ satisfies Condition $(E)$. $\blacksquare$
\begin{cor}
Let $S$ be a pomonoid over which all generators satisfy Condition
$(E)$. Then for any non-empty family $\{A_i~|~i\in I\}$ of right
$S$-posets satisfying Condition $(E), \prod_I A_i$ satisfies
Condition $(E)$.
\end{cor}
\vspace{.2cm} \noindent {\bf Proof.} Using Corollary~\ref{S^I},
$S^I$ satisfies Condition $(E)$ for each non-empty set $I$ and the
result follows by~\cite[Theorem 3.3]{khos:dir}. $\blacksquare$

We remark that  left po-cancellative pomonoids are examples of
pomonoids over which all generators satisfy Condition $(E)$. Since
over such pomonoids, $S_S\times A_S$ satisfies Condition $(E)$, for
every right $S$-poset $A_S$, so the result follows from the former
theorem.
\begin{lem}\label{cycproj}
Let $S$ be a pomonoid and $e^2=e\in S$. Then $\nu(1, e)\subseteq
\ker \lambda_e$ {\rm (}the right monocyclic congruence induced by
the pair $(1, e)${\rm )}.
\end{lem}
\vspace{.2cm} \noindent {\bf Proof.} Suppose $(s, t)\in \nu(1, e)$
and $s\neq t.$ Using Lemma~\ref{induced} there exist \\
$s_1, s_2,
\cdots, s_n, t_1, t_2, \cdots, t_m\in S$ such that
\begin{center}
$s\leq c_1s_1,~d_1s_1\leq c_2s_2,~d_2s_2\leq c_3s_3,~\cdots,~d_ns_n\leq t$;\\
$t\leq p_1t_1,~q_1t_1\leq p_2t_2,~q_2t_2\leq p_3t_3,~\cdots,~q_mt_m\leq s$,
\end{center}
where $(c_i, d_i)\in \{(1, e)\}$, $(p_j, q_j)\in \{(1, e)\}$ for $i=1, 2, \cdots, n$,
$j=1, 2, \cdots, m$. So $ec_i=ed_i$ for each $1\leq i\leq n$ and $ep_j=eq_j$ for each $1\leq j\leq m$ which implies
\begin{center}
$es\leq ec_1s_1=ed_1s_1\leq ec_2s_2=ed_2s_2\leq ~\cdots~\leq ed_ns_n\leq et$;\\
$et\leq ep_1t_1=eq_1t_1\leq ep_2t_2=eq_2t_2\leq ~\cdots~\leq eq_mt_m\leq es$,
\end{center}
therefore $es=et$. $\blacksquare$

We get from Theorem~\ref{mainE} and Lemma~\ref{cycproj} the
following.
\begin{cor}
Let $S$ be a pomonoid. If all generators satisfy Condition $(E)$ then for each $x, y\in S$,
$\nu(1, e)\subseteq \nu(x, y)$ for some idempotent $e\in S$ in case $l(x, y)\neq \emptyset$.
\end{cor}
\section{Pomonoids over which all generators  are  \\ principally weakly po-flat}
The following definition of $A(I)$ for $S$-posets first appeared in~\cite{mahmoudi} and then was used again
in~\cite{Qili2} and~\cite{Qili1}.

Suppose $S$ be a pomonoid and $I$ a proper right ideal of $S$, if $x,y$ and $z$ denote elements not belonging to $S$,
let
\begin{equation}\label{}
   A(I)=(\{x, y\}\times (S-I))\cup (\{z\}\times I)
\end{equation}
and define a right $S$-action on $A(I)$ by
\begin{equation*}
(x, u)s= \left\{
\begin{array}{rl}
(x, us), & \text{if } us\notin I,\\
(z, us), & \text{if } us\in I,
\end{array} \right.
\end{equation*}
\begin{equation*}
(y, u)s= \left\{
\begin{array}{rl}
(y, us), & \text{if } us\notin I,\\
(z, us), & \text{if } us\in I,
\end{array} \right.
\end{equation*}
\begin{center}
$(z, u)s=(z, us)$.
\end{center}
The order on $A(I)$ is defined as follows:
\begin{align*}
(w_1, s)\leq (w_2, t)\Leftrightarrow &(w_1=w_2~\text{and}~s\leq t)~\text{or}\\
                                                   &(w_1\neq w_2, ~s \leq i\leq t~\text{for some}~ i\in I)
\end{align*}
In~\cite{mahmoudi} the authors called its the  \emph{amalgamated coproduct of two copies of S over I} and denoted by
$S\sqcup_I S$.
\begin{lem}\label{A(I)}
Let $S$ be a pomonoid and $I$ be a proper right ideal of $S$. Then
$A(I)$ {\rm (}amalgamated coproduct of two copies of S over I{\rm )}
is a generator in {\bf Pos}-$S$.
\end{lem}
\vspace{.2cm} \noindent {\bf Proof.} Suppose $S$ is a pomonoid, $I$
a right ideal of $S$ and $\iota :I\hookrightarrow S$ is the
inclusion map.
\begin{displaymath}
\xymatrix{
I\ar@{^{(}->}[d]_{\iota} \ar@{^{(}->}[r]^{\iota} & S\ar[d]^{q_1} \ar[rdd]^{\id_S} &\\
S\ar[drr]_{\id_S} \ar[r]^{q_2}   & A(I)\ar@{-->}[dr]^{f}        &\\
  &      & S}
\end{displaymath}
By the universal property of pushout, there exists $S$-homomorphism
$f$ such that $fq_i=\id_S$ for $i=1, 2$. Therefore $S$ is a retract
of $A(I)$ and the result follows by~\cite[Theorem 2.1]{gen}.
$\blacksquare$
\begin{thm}\label{regular}
The following conditions on a pomonoid $S$ are equivalent:
\begin{itemize}
      \item[(i)] All generators in {\bf Pos}-$S$ are principally weakly po-flat;
      \item[(ii)] All right $S$-posets are principally weakly po-flat;
      \item[(iii)] $S$ is a regular pomonoid. That is, for any $s\in S$, there exists $x\in S$ such that $s=sxs$.
\end{itemize}
\end{thm}
\vspace{.2cm} \noindent {\bf Proof.}
$(i) \Rightarrow (iii)$ Let $s\in S$. If $sS=S$, then there exists $x\in S$ such that $sx=1_S$ so $s=sxs$.
Therefore $S$ is regular. If $sS\neq S$ then $I=sS$ is a proper right ideal of $S$.
Then $A(I)$ is a generator by Lemma~\ref{A(I)} and therefore, by hypothesis, principally weakly po-flat.
As $(x, 1_S)s\leq (y, 1_S)s$ in $A(I)$ we get $(x, 1_S)\otimes s\leq (y, 1_S)\otimes s$ in $A(I)\otimes {}_{S} S$ and,
by hypothesis, just as well in $A(I)\otimes {}_{S} Ss$, i.e. we have an $S$-tossing in $A(I)\otimes {}_{S} Ss$
connecting $((x, 1_S), s)$ and $((y, 1_S), s)$ by~\cite[Lemma 1.2]{cyclic}.
By considering $S$-tossing and the $S$-action and the order of $A(I)$, it can be followed $S$ is regular.

$(iii) \Rightarrow (ii)$ Let $S$ be regular. Take $A_S\in$ {\bf
Pos}-$S$ and let $as\leq a's$ for $a, a'\in A_S$, $s\in S$. Since
$S$ is regular there exists $x\in S$ with $s=sxs$. Then
$$a\otimes s\leq a\otimes sxs\leq as\otimes xs\leq a's\otimes xs\leq a'\otimes sxs\leq a'\otimes s$$
in $A\otimes {}_{S} Ss$. This proves that $A_s$ is principally weakly po-flat by~\cite[Proposition 2.7]{cyclic}.

$(ii) \Rightarrow (i)$ This is trivial. $\blacksquare$
\section{Pomonoids over which all generators  are weakly po-flat}
Our next aim is characterizing pomonoids over which all generators
are weakly po-flat. We record the following theorem.
\begin{thm}~\cite[Theorem 3.12]{sflat poflat}\label{W}
A right $S$-poset $A$ is weakly po-flat if and only if it is principally weakly po-flat and
satisfies the following Condition:

$(W)$ If $ax\leq a'y$ for $a, a'\in A$, $x, y\in S$, then there exist $a''\in A$, $p\in Sx$, $q\in Sy$,
such that $p\leq q$, $ax\leq a''p$ and $a''q\leq a'y$.
\end{thm}
Similar to the definition in S-acts (see~\cite[Definition 3.2]{sed:gen}), we have
\begin{defn}\label{almost}
A right $S$-poset $A_S$ is called \emph{almost weakly po-flat} if $A_S$ is principally weakly po-flat
and satisfies condition

$(W')$ if $ax\leq a'y$ and $Sx\cap (Sy]\neq \emptyset$ for $a, a'\in A_S$, $x, y\in S$,
then there exist $a''\in A_S$, $p\in Sx$ and $q\in Sy$ such that $p\leq q$, $ax\leq a''p$ and $a''q\leq a'y$.
\end{defn}

Obviously, Condition $(W')$ follows from Condition $(W)$. So we have
the implications:\\ weakly po-flat $\Rightarrow$  almost weakly
po-flat $\Rightarrow$ principally weakly po-flat.

Note that if $S$ is weakly right reversible, the concepts almost weak po-flatness and weak po-flatness coincide.\\
In the next lemma we show that the property weak po-flatness is preserved under retractions.
\begin{lem}\label{retwpoflat}
Let $S$ be a pomonoid. Then any retract of a weakly po-flat $S$-poset is weakly po-flat.
\end{lem}
\vspace{.2cm} \noindent {\bf Proof.} Assume $B_S$ is a retract of a
right $S$-poset $A_S, \xymatrix{A\ar@<0.5ex>[r]^{f}&
B\ar@<0.5ex>[l]^{g}}, fg=\id_{B_S}, A_S$ is weakly po-flat and
$bs\leq b't$ for $b, b'\in B_S, s, t\in S$. Then $g(b)s\leq g(b')t$
in $A_S$, using the definition of weakly po-flatness of $A_S$ we
have $g(b)\otimes s\leq g(b')\otimes t$ in $A\otimes {}_{S} (Ss\cup
St)$. So there exist $a_1, a_2, \cdots, a_n\in A$, $s_1, s_2,
\cdots, s_n,t_1, t_2, \cdots, t_n\in S$ and $u_2, u_3, \cdots,
u_n\in Ss\cup St$ such that
\begin{center}
$g(b)\leq a_1s_1$~~~~~~~~~~~~~~~~~~~~\\
$a_1t_1\leq a_2s_2\qquad s_1s\leq t_1u_2$\\
$a_2t_2\leq a_3s_3\qquad s_2u_2\leq t_2u_3$\\
$\cdots~~~~\qquad ~~~~~\cdots$\\
$a_nt_n\leq g(b')\qquad s_nu_n\leq t_nt$
\end{center}
Using that $f$ is $S$-homomorphism and applying $f$, we get
\begin{center}
$b=fg(b)\leq f(a_1)s_1$~~~~~~~~~~~~~~~~~~~~\\
$f(a_1)t_1\leq f(a_2)s_2\qquad s_1s\leq t_1u_2$\\
$f(a_2)t_2\leq f(a_3)s_3\qquad s_2u_2\leq t_2u_3$\\
$~~~\cdots~~\qquad ~~~~~~~~~~~~~~\cdots$\\
$f(a_n)t_n\leq fg(b')=b'\quad s_nu_n\leq t_nt$
\end{center}
Now~\cite[Lemma 1.2]{cyclic} implies $b\otimes s\leq b'\otimes t$ in
$B\otimes {}_{S} (Ss\cup St)$. Thus $B_S$ is weakly po-flat.
$\blacksquare$
\begin{thm}\label{wpoflat}
The following assertions are equivalent for a pomonoid $S$:
\begin{itemize}
      \item[(i)] all generators are weakly po-flat;
      \item[(ii)] $S_S\times A_S$ is weakly po-flat for each right $S$-poset $A_S$;
      \item[(iii)] a right $S$-poset $A_S$ is weakly po-flat if $\Pos_S(A_S,S_S)\neq \emptyset$;
      \item[(iv)] all right $S$-posets are almost weakly po-flat.
\end{itemize}
\end{thm}
\vspace{.2cm} \noindent {\bf Proof.}
Assertions $(i)$, $(ii)$ and $(iii)$ are equivalent by Theorem~\ref{gen} and Lemma~\ref{retwpoflat}.

$(i) \Rightarrow (iv)$ Let $A_S$ be a right $S$-poset. By assumption
all generators are weakly po-flat, so they are principally weakly
po-flat which Theorem~\ref{regular} implies $S$ is regular and then
$A_S$ is principally weakly po-flat. Next, suppose that $as\leq a't$
and $Ss\cap (St]\neq \emptyset$ for $a, a'\in A_S$, $s, t\in S$.
This makes $xs\leq yt$ for some $x, y\in S$. So $(x, a)s\leq (y,
a')t$ in $S_S\times A_S$. Since $S_S\times A_S$ is generator, then
it is weakly po-flat by assumption. Thereby $S_S\times A_S$
satisfies Condition $(W)$, which makes $A_S$ satisfies condition
$(W')$. Therefore $A_S$ is almost weakly po-flat.

$(iv) \Rightarrow (i)$ Let $A_S$ be a generator and $f:A_S\arrow
S_S$ be an epimorphism. Since by assumption $A_S$ is almost weakly
po-flat, by Theorem~\ref{W} it suffices to show that $A_S$ satisfies
Condition $(W)$. Suppose that $as\leq a't$ for $a, a'\in A_S$, $s,
t\in S$. Since $f(a)s\leq f(a')t$ and $f(a), f(a')\in S$ we get
$Ss\cap (St]\neq \emptyset$ so  Condition $(W')$ is satisfied,
thereby there exist $a''\in A_S$, $p\in Ss$ and $q\in St$ such that
$as\leq a''p$, $a''q\leq a't$ and $p\leq q$. Then $A_S$ is weakly
po-flat. $\blacksquare$

The proof of the next theorem exactly follows from the corresponding proof for $S$-acts
(see~\cite[Corollary 3.5]{character}).
\begin{thm}
Let $S$ be a pomonoid with a left zero. Then all generators are weakly po-flat if and only if all right
$S$-posets are weakly po-flat.
\end{thm}

In the next corollary, using Theorem~\ref{wpoflat}, we reduce this condition on $S$ to weakly right reversibility.
\begin{cor}
Let $S$ be a weakly right reversible pomonoid. Then all generators are weakly po-flat if and only if all right
$S$-posets are weakly po-flat.
\end{cor}
\vspace{.2cm} \noindent {\bf Proof.} It follows from
Theorem~\ref{wpoflat} and the remark after Definition~\ref{almost}.
$\blacksquare$

In 1960, left (right) coherent rings were characterized by Chase~\cite{Chase} as rings over which
every direct product of flat right (left) modules is flat. We remind the reader that a monoid $S$ is
\emph{(weakly) left coherent} if products of (weakly) flat right $S$-acts are (weakly) flat.
Weakly left coherent monoids are characterized in~\cite[Theorem 2.11]{seda:coher}.
It is worth noting that from another aspect (weakly) left coherent monoids are defined in~\cite{Gould}
as monoids over which every finitely generated (left ideal) subact of ($S$) every finitely presented left $S$-act
is finitely presented.

Let $S$ be a pomonoid. It is easily seen that coproducts of families of (weakly) po-flat right $S$-posets are again
(weakly) po-flat. Taking our lead from ring theory, we shall call $S$ \emph{weakly left coherent}
if all direct products of non-empty families of weakly po-flat right $S$-posets are weakly po-flat.

Now we are going to show that pomonoids over which all generators are weakly po-flat are regular
and weakly left coherent. Furthermore we characterize weakly left coherent pomonoids in special case.
In what follows we collect the results that are needed to reach this goal.

Recall that a finitely generated left $S$-poset ${}_{S}B$ is called \emph{finitely definable $(FD)$}
if the $S$-homomorphism $S^I\otimes {}_{S} B\arrow B^I$, $(s_i)\otimes b\rightsquigarrow (s_ib)$,
is order-embedding for all nonempty sets $I$. For each $(p, q)\in D(S)$,
$$\{(u, v)\in D(S)~|~ \exists w\in S~; u\leq wp, wq\leq v\}$$
is a left $S$-poset and will be denoted by $\widehat{S(p, q)}$ from now on. Clearly $\widehat{S(p, q)}$
contains the cyclic $S$-poset S$(p, q)$. Moreover, if $Ss\cap (St] ̸\neq\emptyset$,
\begin{equation}\label{}
  \{(as, a't)\in Ss\times St~|~ as\leq a't\}
\end{equation}
is denoted by $H(s, t)$ (see~\cite{khos:dir}).
\begin{thm}\cite[Theorem 2.7]{khos:dir}\label{khos2.7}
The following are equivalent for a pomonoid $S$:
\begin{itemize}
      \item[(i)] $S^I$ is weakly po-flat right $S$-poset for each non-empty set $I$;
      \item[(ii)] every finitely generated left ideal of $S$ is $FD$;
      \item[(iii)] $Sx$ is $FD$ for each $x\in S$, and\\
for every $x, y\in S$, if $Sx\cap (Sy]\neq \emptyset$, then $H(x, y)\subseteq \widehat {S(p, q)}$
for some $(p, q)\in H(x, y)$.
\end{itemize}
\end{thm}
\begin{lem}\cite[Proposition 4.8]{Inde}\label{PP}
An ordered monoid $S$ is a left $PP$ pomonoid if and only if for every $a\in S$ there exists an idempotent $e$ of $S$
such that $a=ea$ and $sa\leq ta$ implies $se\leq te$ for $s, t\in S$.
\end{lem}
\begin{thm}\cite[Corollary 3.17]{sflat poflat}\label{PP wpo-flat}
Let $S$ be a left $PP$ pomonoid and $A_S$ a right $S$-poset. Then the following conditions are equivalent:
\begin{itemize}
      \item[(i)] $A_S$ is weakly po-flat;
      \item[(ii)] For any $a, a'\in A_S$, $x, y\in S$, if $ax\leq a'y$, then there exist $a''\in A_S$, $x_1, y_1\in S$
      and $u, v\in E(S)$ such that
$$ux=x,~~vy=y,~~x_1x\leq y_1y$$
$$au\leq a''x_1,~~a''y_1\leq a'v.$$
\end{itemize}
\end{thm}
In the next theorem weakly left coherent pomonoids are characterized in special case.
\begin{thm}\label{wlcohere}
Let $S$ be a left $PP$ pomonoid. Then the following statements are equivalent:
\begin{itemize}
      \item[(i)] for every non-empty set $I$, $S^I$ is a weakly po-flat right $S$-poset;
      \item[(ii)] every finitely generated left ideal of $S$ is $FD$;
      \item[(iii)] for every $x, y\in S$, if $Sx\cap (Sy]\neq \emptyset$, then
      $H(x, y)\subseteq \widehat{S(p, q)}$ for some $(p, q)\in H(x, y)$.
      \item[(iv)] $S$ is weakly left coherent.
\end{itemize}
\end{thm}
\vspace{.2cm} \noindent {\bf Proof.}
Implications $(i)\Rightarrow (ii) \Rightarrow (iii)$ follow from Theorem~\ref{khos2.7} and
$(iv)\Rightarrow (i)$ is trivial.

$(iii)\Rightarrow (iv)$ Let $\{A_i~|~i\in I\}$ be a non-empty family of weakly po-flat right $S$-posets.
Using Theorem~\ref{PP wpo-flat} we show $A=\prod_{i\in I} A_i$, is weakly po-flat. Suppose $x, y\in S$,
$(a_i), (a'_i)\in A$ and $(a_i)x\leq (a'_i)y$. Then for each, $i\in I$ $a_ix\leq a'_iy$ in $A_i$ and
so by Theorem~\ref{PP wpo-flat} there exist $\hat{a_i}\in A_i$, $x_i, y_i\in S$ and $u_x, u_y\in E(S)$ such that
$$x=u_xx,~~y=u_yy,~~x_ix\leq y_iy$$
$$a_iu_x\leq \hat{a_i}x_i,~~\hat{a_i}y_i\leq a'_iu_y$$
Since $x_ix\leq y_iy$, then $Sx\cap (Sy]\neq \emptyset$, by
assumption $H(x, y)\subseteq \widehat{S(p, q)}$ for some $(p, q)\in
H(x, y)$. Therefore $(x_ix, y_iy)\in \widehat{S(p, q)}$  for some
$(p, q)\in H(x, y)$, so there exists $z_i\in S$ such that $x_ix\leq
z_ip$ and $z_iq\leq y_iy$. On the other hand, since $(p, q)\in H(x,
y)$, then there exist $c, d\in S$ such that $p=cx$, $q=dy$ and
$cx\leq dy$. By Lemma~\ref{PP}  and the proof of Theorem~\ref{PP
wpo-flat} we know that $sx\leq tx$ implies $su_x\leq tu_x$ for $s,
t\in S$, similarly $y$ has this property. So for each $i\in I$ we
can conclude, from the inequalities $x_ix\leq z_icx$ and $z_idy\leq
y_iy$, that $x_iu_x\leq z_icu_x$ and $z_idu_y\leq y_iu_y$. Define
$a''_i=\hat{a_i}z_i$ for each $i\in I$, $x_1=cu_x$ and $y_1=du_y$.
Then
$$(a''_i)x_1=(\hat{a_i}z_i)cu_x=(\hat{a_i}z_icu_x)\geq (\hat{a_i}x_iu_x)\geq (a_iu_xu_x)=(a_i)u_x.$$
So $(a_i)u_x\leq (a''_i)x_1$. Similarly it follows that
$(a''_i)y_1\leq (a'_i)u_y$. Finally, $x_1x=cu_xx=cx\leq
dy=du_yy=y_1y$. The proof is now complete. $\blacksquare$
\begin{cor}
Let $S$ be a pomonoid over which all generators are weakly po-flat.
Then $S$ is a regular and weakly left coherent pomonoid or equivalently $S$ is regular and for every $x, y\in S$,
if $Sx\cap(Sy]\neq \emptyset$, then $H(x, y)\subseteq \widehat{S(p, q)}$ for some $(p, q)\in H(x, y)$.
\end{cor}
\vspace{.2cm} \noindent {\bf Proof.} Since all generators are
principally weakly po-flat, all right $S$-posets are principally
weakly po-flat and we conclude $S$ is a regular pomonoid (see
Theorem~\ref{regular}). By Theorem~\ref{wlcohere} a left $PP$
pomonoid $S$ is weakly left coherent if and only if $S^I$ is weakly
po-flat for each non-empty set $I$ or equivalently for every $x,
y\in S$, if $Sx\cap(Sy]\neq \emptyset$, then $H(x, y)\subseteq
\widehat{S(p, q)}$ for some $(p, q)\in H(x, y)$. Now, using
Corollary~\ref{S^I}, $S^I$ is weakly po-flat for each nonempty set
$I$. Therefore, using the fact that regular pomonoids are left $PP$,
$S$ is a regular and weakly left coherent pomonoid. $\blacksquare$
\section{Pomonoids over which all generators are $I$-regular}
In this section, the concept of an $I$-regular $S$-poset is reminded and characterization of pomonoids over which
all generators are $I$-regular is given. Similar to regular acts, the concept of $I$-regular $S$-poset
was introduced in~\cite{Inde}. Let $A_S$ be an $S$-poset. An element $a\in A_S$ is called \emph{$I$-regular}
if there exists an $S$-homomorphism $f:aS\arrow S$ such that $af (a)=a$. An $S$-poset $A_S$ is called \emph{$I$-regular}
if all elements of $A_S$ are $I$-regular.
\begin{lem}\label{retIreg}
Let $S$ be a pomonoid. Then any retract of an $I$-regular $S$-poset
is $I$-regular.
\end{lem}
\vspace{.2cm} \noindent {\bf Proof.}
Assume $B_S$ is a retract of a right $S$-poset $A_S,
 \xymatrix{A\ar@<0.5ex>[r]^{\gamma}&
B\ar@<0.5ex>[l]^{\pi}}, \pi \gamma=\id_{B_S}$ and $A_S$ is
$I$-regular.
 We must show that $B_S$ is $I$-regular. Let $b\in B_S$, using~\cite[Theorem 4.2 and Corollary 4.3]{Inde},
 there exists an $S$-homomorphism $\beta:\gamma(b)S\arrow eS$ such that $\beta (\gamma(b))=e$, $\gamma (b)=\gamma (b)e$
 where $e\in E(S)$. Now we define
\begin{align*}
\alpha &: \pi \gamma(b)S\arrow S \\
         &\qquad b\mapsto \beta(\gamma(b))
\end{align*}
it is clear that $\alpha$ is an $S$-homomorphism and $\alpha(b)=\beta(\gamma(b))=e\in E(S)$.
It suffices to show that $b\beta(\gamma(b))=b$.
$$b\beta(\gamma(b))=\pi \gamma(b)\beta(\gamma(b))=\pi \gamma(b)e=\pi \gamma(b)=b.$$
Therefore $B_S$ is $I$-regular. $\blacksquare$
\begin{thm}
Let $S$ be a pomonoid. The following assertions are equivalent.
\begin{itemize}
      \item[(i)] all generators are $I$-regular;
      \item[(ii)] $S_S\times A_S$ is $I$-regular for every right $S$-poset $A_S$;
      \item[(iii)] a right $S$-poset $A_S$ is $I$-regular if $\Pos_S(A_S,S_S)\neq \emptyset$
      \item[(iv)] for each right subannihilator congruence $\rho$, $S/\rho$ is $I$-regular.
\end{itemize}
\end{thm}
\vspace{.2cm} \noindent {\bf Proof.} Assertions $(i)$, $(ii)$ and
$(iii)$ are equivalent by Theorem~\ref{gen} and Lemma~\ref{retIreg}.
Implications $(i) \Rightarrow (iv)$ and $(iv) \Rightarrow (i)$ is
similar to the act case (see~\cite{sed:gen}). $\blacksquare$

Recall that a pomonoid $S$ is called a \emph{right $PP$ pomonoid} if the ideal $xS$ is projective for all $x\in S$.
Let $A_S$ be a right $S$-poset over a pomonoid $S$. If we denote by $\Con(A_S)$ the set of all congruences on the
$S$-poset $A_S$, then clearly $(\Con(A_S),  \subseteq, \cap)$ is a commutative band with the identity element
$A_S\times A_S$.
Our next aim is to show that when every product of $I$-regular $S$-posets is $I$-regular.
To reach this goal we need some results that being entirely similar to the $S$-act case,
which can be found in~\cite{seda:principal}.
\begin{thm}\label{reg1}
Let $S$ be a pomonoid. The diagonal $S$-poset $D(S)$ is $I$-regular if and only if
\begin{itemize}
      \item[(a)] $S$ is a right PP pomonoid and
      \item[(b)] the set $R=\{\ker \lambda_e~|~e\in E(S)\}\cup \{S\times S\}$  is an subpomonoid of
      $T=(\Con(S_S),\subseteq, \cap)$.
\end{itemize}
\end{thm}
\vspace{.2cm} \noindent {\bf Proof.} The proof is similar to the
$S$-act case. $\blacksquare$
\begin{thm}\label{reg2}
The following are equivalent for a right $PP$ pomonoid $S$:
\begin{itemize}
     \item[(i)] every finite product of $I$-regular $S$-posets is $I$-regular;
     \item[(ii)] $S^n$ is $I$-regular for every $n\in \mathbb{N}$;
     \item[(iii)] the diagonal $S$-poset $D(S)$ is $I$-regular.
\end{itemize}
\end{thm}
\vspace{.2cm} \noindent {\bf Proof.} Using Theorem~\ref{reg1} and
similar to the $S$-act case. $\blacksquare$
\begin{thm}\label{prod reg}
The following are equivalent for a right $PP$ pomonoid $S$:
\begin{itemize}
    \item[(i)] every product of $I$-regular $S$-posets is $I$-regular;
    \item[(ii)] $S^I$ is $I$-regular for each nonempty set $I$;
\item[(iii)] for each family $\{e_i~|~i\in I\}$ of idempotents of $S$, there exists an idempotent $e\in S$ such that
$\bigcap_{i\in I}\ker \lambda_{e_i}=\ker \lambda_e$
\end{itemize}
\end{thm}
\vspace{.2cm} \noindent {\bf Proof.} Using Theorem~\ref{reg2} and
similar to the $S$-act case. $\blacksquare$
\begin{cor}
Let $S$ be a pomonoid over which all generators are $I$-regular. Then every product of $I$-regular
$S$-posets is $I$-regular.
\end{cor}
\vspace{.2cm} \noindent {\bf Proof.} Since $S$ is a generator then
it is an $I$-regular $S$-poset. Hence $S$ is a right $PP$ pomonoid.
On the other hand, using Corollary~\ref{S^I} and
Lemma~\ref{retIreg}, $S^I$ is $I$-regular for each nonempty set $I$
and by Theorem~\ref{prod reg} every product of $I$-regular
$S$-posets is $I$-regular. $\blacksquare$
\section{Pomonoids over which all generators are po-torsion free}
In~\cite{sflat poflat} and~\cite{bulman flatness} the obvious definitions of torsion freeness
and po-torsion freeness are presented and, by analogy with the situation for $S$-acts,
it is shown that principally weakly flatness and principally weakly po-flatness strictly imply, respectively,
torsion freeness and po-torsion freeness, for a given $S$-poset. Furthermore,
it is shown that po-torsion freeness and torsion freeness are incomparable properties.

Now we are going to be characterized pomonoids over which all generators are po-torsion free.
An element $c$ of a pomonoid $S$ will be called \emph{right po-cancellable} if for all $s, s'\in S$, $sc\leq s'c$
implies $s\leq s'$. An $S$-poset $A_S$ will be called \emph{po-torsion free} if $ac\leq a'c$ implies $a\leq a'$
whenever $a, a'\in A_S$ and $c$ is a right po-cancellable element of $S$. Right po-cancellable elements are
right cancellable in the usual sense, but the converse is not true.
\begin{thm}
The following conditions on a pomonoid $S$ are equivalent:
\begin{itemize}
   \item[(i)] all generators in {\bf Pos}-$S$ are po-torsion free;
   \item[(ii)] all right $S$-posets are po-torsion free;
   \item[(iii)] every right po-cancellable element of $S$ is right invertible.
\end{itemize}
\end{thm}
\vspace{.2cm} \noindent {\bf Proof.}
Implications $(iii) \Rightarrow (ii)$ and $(ii) \Rightarrow (i)$ are trivial.

$(i) \Rightarrow (iii)$ Suppose that $s\in S$ is right
po-cancellable but not right invertible. Then $I=sS$ is a proper
right ideal of $S$ and thus $A(I)$ is a generator by
Lemma~\ref{A(I)}, which, by hypothesis, is po-torsion free. But $(x,
1_S)s\leq (y, 1_S)s$ which implies $(x, 1_S)\leq (y, 1_S)$. So there
exists $st\in I=sS$, such that $1_S\leq st\leq 1_S$. Therefore
$st=1_S$ which is a contradiction. $\blacksquare$

\end{document}